\documentclass[12pt]{article}

\usepackage[T1]{fontenc}
\usepackage{avant}
\usepackage[cp1250]{inputenc}
\usepackage{amsmath,amssymb,amsthm}
\usepackage{amscd}
\usepackage{etoolbox}
\usepackage{lipsum}
\usepackage{tikz}
\usepackage{tabularx}

\usetikzlibrary{matrix}

\newtheorem{theorem}{Theorem}

\newtheorem{lemma}{Lemma}
\newtheorem{definition}{Definition}

\newtheorem{remark}{Remark}

\title{On the Hirzebruch-Kobayashi-Ono proportionality principle and the non-existence of compact solvable Clifford-Klein forms of certain homogeneous spaces}

\author{Maciej Boche\'nski and Aleksy Tralle\\
 Department of Mathematics and Computer Science,\\
  University of Warmia and Mazury, \\
  S\l\/oneczna 54, 10-710, Olsztyn, Poland\\
  E-mail: mabo@matman.uwm.edu.pl (MB), \\ 
  aleksymath46@gmail.com \, (AT)}

\begin{document}
\maketitle{}
\begin{abstract}This article continues a line of research aimed at solving an important problem of T. Kobayashi of the existence of compact Clifford-Klein forms of reductive homogeneous spaces. We contribute to this topic by showing that almost all symmetric spaces and 3-symmetric spaces do not admit solvable compact Clifford-Klein forms (with  several possible exceptions). Our basic tool is a combination of the Hirzebruch-Kobayashi-Ono proportionality principle with the theory of syndetic hulls. Using this, we prove a general theorem which yields a sufficient condition for the non-existence of compact solvable Clifford-Klein forms.
\vskip6pt

\end{abstract}
\section{Introduction}\label{sec:intro}
Assume that we are given a non-compact homogeneous space $G/H$ of a reductive real Lie group $G$ and a closed subgroup $H\subset G$. Suppose that there exists a discrete subgroup $\Gamma\subset G$ which acts properly,  co-compactly and freely on $G/H$ by left translations. The quotient $\Gamma\backslash G/H$ is called a compact Clifford-Klein form. The problem of determining which reductive homogeneous spaces admit such forms goes back to Calabi and Markus \cite{CM62}  and was formulated as a research program by T. Kobayashi \cite{ko-ma}. The main inspiration for the research presented in this paper are the works \cite{ko-ma}, \cite{ko}, \cite{kob92}, \cite{kkb} and \cite{koy}. In \cite{kob92}  the question was explicitly formulated as follows: {\it when does there exist a compact Clifford-Klein form of homogeneous space of reductive type?}
In general, the problem is still not solved, neither for symmetric spaces, nor for other wide classes of homogeneous spaces. In  \cite{koy} various  ways to attack it are discussed (for example, Theorems 2.2.1, 2.3.1, 2.4.1, Corollary 3.5.9, Corollaries 3.6.4 and 3.6.5).  
Also, \cite{koy} contains  several important conjectures on the (non)existence of compact Clifford-Klein forms of symmetric spaces, as well as a description of general methods of the theory of Clifford-Klein forms: criteria of properness of the group actions  and obstructions. In this article we focus on a question: {\it what algebraic conditions should be satisfied by a discrete subgroup $\Gamma\subset G$ if it yields a compact Clifford-Klein form?} This question is partially motivated by a theorem of Y. Benoist \cite{b}. It was shown in \cite{b} that no nilpotent group can act properly and co-compactly on a non-compact homogeneous space $G/H$ with semisimple $G$. Note that not much is known about the possible algebraic properties of discrete groups that may yield compact Clifford-Klein forms: the main examples are lattices in closed subgroups $L$ which are reductive in a given reductive Lie group $G$ and their deformations \cite{kass},\cite{ko-ma}. In our line of thinking the  results of the present article continue  \cite{bt}. In \cite{bt} we have proved a generalization of Benoist's theorem to solvable $\Gamma$ under the assumption that $G$ is reductive and  $H$ is a semisimple part of a centralizer of a torus in $G$. It is worth noting that our approach to the problem differs from that of \cite{bt} and might be of independent interest. We use a result of T. Kobayashi and K. Ono \cite{ko} (which we call the Hirzebruch-Kobayashi-Ono proportionality principle), together with an approach to the Clifford-Klein forms via syndetic hulls \cite{w-a}. Thus, we combine topological methods with methods of the theory of algebraic groups, in order to understand the algebraic properties of $\Gamma$. This combination of methods yields new series of examples of homogeneous spaces $G/H$ with no compact {\it solvable} Clifford-Klein forms: these are non-compact absolutely simple symmetric spaces (with four possible exceptions) and non-compact absolutely simple 3-symmetric spaces (in the sense of Wolf and Gray \cite{wg1}, \cite{wg2}), with one possible exception.  Note that we define  a solvable compact Clifford-Klein form $\Gamma\backslash G/H$ as a quotient of $G/H$ by a discrete solvable subgroup $\Gamma\subset G$ acting properly, freely  and co-compactly on $G/H$ by left translations.

Let $K$ be a maximal compact subgroup of $G$. Choose it in a way that $K_H=K\cap H$ is a maximal compact subgroup of $H$. Consider the corresponding Lie algebras $\mathfrak{g},\mathfrak{k},\mathfrak{h}$ and $\mathfrak{k}_H$. Take a Cartan decomposition $\mathfrak{g}=\mathfrak{k}+\mathfrak{p}$. There is an orthogonal (with respect to the Killing form of $\mathfrak{g}$) decomposition
$$\mathfrak{g}=\mathfrak{k}_H+\mathfrak{k}_H^{\perp}+\mathfrak{p}_H+\mathfrak{p}_H^{\perp}$$
where $\mathfrak{p}_H=\mathfrak{p}\cap\mathfrak{h}.$ Also denote by $H^*(\mathfrak{g},\mathbb{R})$ the Lie algebra cohomology of $\mathfrak{g}$ (\cite{K1}). 

\begin{theorem}
Let $G/H$ be a non-compact  homogeneous space of a semisimple linear Lie group $G.$ Assume that $G/H$ is of reductive type. Then $G/H$ admits a compact solvable Clifford-Klein form only if the following conditions are all satisfied:

\begin{enumerate}
  \item $H^n(\mathfrak{g},\mathbb{R})\not=0$ for $n=\dim\mathfrak{k}+\dim\mathfrak{p}_H$.
	\item $\operatorname{rank}\,G>\operatorname{rank}\,H$.
\end{enumerate}
\label{twgg}
\end{theorem}

As an application of this theorem, we get the non-existence of solvable compact Clifford-Klein forms for almost all simple symmetric spaces. 

\begin{theorem}\label{cwgg} Let $G/H$ be a symmetric space of an absolutely simple non-compact connected linear real Lie group $G$, and $(\mathfrak{g},\mathfrak{h})$ be the corresponding symmetric pair.  If the pair $(\mathfrak{g},\mathfrak{h})$ is different from
$$(\mathfrak{so}(n,m),\mathfrak{so}(n-k,m)\oplus \mathfrak{so}(k)),\,\,\text{with}\,\,0<n\leq m;\,m,\, n \,\text{even}, \ k \ \text{odd},$$
$$(\mathfrak{su}(2p,2q),\mathfrak{sp}(p,q)),\,(\mathfrak{su}(2m-1,2m-1),\mathfrak{so}^*(4m-2)), (\mathfrak{e}_{6(2)}, \mathfrak{sp}(3,1)),$$
then $G/H$ does not admit a solvable compact Clifford-Klein form.
\end{theorem}

In the same way, we get the following.
\begin{theorem}\label{thm:3-sym} Let $G/H$ be  a 3-symmetric space of an absolutely simple non-compact connected linear real Lie group $G$, and $(\mathfrak{g},\mathfrak{h})$ be the corresponding 3-symmetric pair. If the pair $(\mathfrak{g},\mathfrak{h})$ is different from
 $$(\mathfrak{so}(3,5),\mathfrak{g}_{2(2)}),$$
then $G/H$ does not admit a solvable compact Clifford-Klein form.
 \end{theorem}
\begin{remark}{\rm Homogeneous spaces $G/H$ of complex simple Lie groups $G$ are not covered by Theorems \ref{cwgg} and \ref{thm:3-sym}}.
\end{remark}
Let us mention topological obstructions to the existence of compact Clifford-Klein forms found by Y. Morita \cite{mor2}, \cite{mor1} and N. Tholozan \cite{Th} (applicable to the general Kobayashi's problem,  not only to the solvable compact Clifford-Klein forms). Note that the Hirzebruch-Kobayashi-Ono proportionality principle has already been used  \cite{ko-ma}, Proposition 4.10, and \cite{ko}, Corollary 5. 

Finally, let us explain the following. Recall that a discrete group $\Gamma$ is {\it amenable} if it admits a finitely additive probability measure. In the context of this article, the term ``amenable'' is almost equivalent to ``solvable''. Indeed, 
 by the Tits alternative if there are no compact solvable Clifford-Klein forms in the considered cases, there are also no amenable ones (because the free group on two generators is non-amenable).

   \vskip6pt
   \noindent {\bf Acknowledgment.} We thank Dmitri Alekseevsky, Ioannis Chrysikos and Jarek K\c edra for answering our questions. We are grateful to the anonymous referee for careful reading of the manuscript and many valuable suggestions which essentially improved the presentation of the results. Our special thanks go to him/her for pointing out the necessity of Lemma 2 and pointing out a mistake in the first version of the manuscript.  
	The first named author acknowledges the support of the National Science Center (grant NCN no. 2018/31/D/ST1/00083). The second named author was supported by NCN grant  2018/31/B/ST1/00053.

\section{Preliminaries}\label{sec:prelim}\subsection{Homogeneous spaces and their Clifford-Klein forms}\label{subsec:c-k}
Throughout this paper we use the basics of Lie theory without further explanations. One can consult \cite{K}. We denote Lie groups by $G, H,...$, and their Lie algebras by the corresponding Gothic letters $\mathfrak{g},\mathfrak{h}...$.  The symbol $\mathfrak{g}^{\mathbb{C}}$ denotes the complexification of a real Lie algebra $\mathfrak{g}$. We also use relations between Lie groups and algebraic groups following \cite{hum} and \cite{w-book}. If ${\bf G}\subset GL(n,\mathbb{C})$ is  an algebraic $\mathbb{R}$-group, then $G={\bf G}_{\mathbb{R}}={\bf G}\cap GL(n,\mathbb{R})$ is a Lie group with a finite number of connected components. If $L$ is a subgroup of $G$ then by $\bar{L}$ we denote the Zariski closure of $L$.

In this article we use various results dealing with  homogeneous spaces $G/H$ of {\it reductive type}. Recall this notion following \cite{ko-ma}. By a real reductive linear Lie group we mean a real linear Lie group $G$ contained in a connected complex reductive Lie group $G^{\mathbb{C}}$ whose Lie algebra $\mathfrak{g}^{\mathbb{C}}$ is isomorphic to $\mathfrak{g}\otimes_{\mathbb{R}}\mathbb{C}$.  We use without further explanations the notions of the {\it Cartan involution} $\theta$ and the {\it Cartan decomposition} $\mathfrak{g}=\mathfrak{k}+\mathfrak{p}$. Let $H$ be a closed subgroup of a real reductive linear group $G$ and $K$ denote the maximal compact subgroup in $G$ corresponding to $\mathfrak{k}$. We say that $H$ is reductive in $G$ and $G/H$ is of {\it reductive type} if there exists a Cartan involution $\theta$ of $G$ such that $G$ has a polar decomposition $H=(H\cap K)\cdot\exp(\mathfrak{h}\cap\mathfrak{p})$ and the connected subgroup $H^{\mathbb{C}}\subset G^{\mathbb{C}}$ corresponding to $\mathfrak{h}^{\mathbb{C}}$ is closed. Note that semisimple linear real Lie groups are reductive. For simplicity, in this article we assume that $G$ is semisimple and linear.
 
 One has the  {\it Iwasawa decomposition} 
$$\mathfrak{g}=\mathfrak{k}+\mathfrak{a}+\mathfrak{n},$$ 
and, on the Lie group level, the {\it global Iwasawa decomposition} $G=\mathit{KAN}$ which is a topological decomposition into a direct product.
Using the Cartan involution, one can write the {\it compatible Cartan decomposition}
  $$\mathfrak{h}=\mathfrak{k}_H+\mathfrak{p}_H, \ \mathfrak{k}_H=\mathfrak{k}\cap\mathfrak{h}, \ \mathfrak{p}_H=\mathfrak{p}\cap\mathfrak{h},$$
  on the Lie algebra level and one can choose a maximal abelian subalgebra $\mathfrak{a}$ in $\mathfrak{p}$ so that $\mathfrak{a}_{H}:=\mathfrak{a}\cap \mathfrak{p}_{H}$ is a maximal abelian subalgebra in $\mathfrak{p}_{H}$.

  Let $X$ be a Hausdorff topological space and $\Gamma$ a topological group acting on $X$. We say that an action of $\Gamma$ on $X$ is {\it proper} if for any compact subset $S\subset X$ the set
  $$\{\gamma\in\Gamma\,|\,\gamma(S)\cap S\not=\emptyset\}$$
  is compact. In particular, this article is devoted to the proper actions of a discrete subgroup $\Gamma\subset G$ on $G/H$ by left translations.  
  
The main result of this article deals with {\it semisimple symmetric spaces}. They are homogeneous spaces $G/H$, where $G$ is a semisimple linear Lie group, and $H$ is a closed subgroup in $G$ such that
	\begin{equation}
  G^{\sigma}_0\subset H\subset G^{\sigma}. 
	\label{equo3}
	\end{equation}
  Here $G^{\sigma}=\{g\in G\,|\,\sigma(g)=g\}$ is the subgroup of fixed points of an involutive automorphism $\sigma$ of $G$, and $G^{\sigma}_0$ denotes the identity connected component of $G^{\sigma}$. Note that we don't assume that $G/H$ is Riemannian (this case is well known: Riemannian symmetric spaces always have compact Clifford-Klein forms by a classical result of Borel \cite{Bor63}). Any symmetric space $G/H$ determines a {\it symmetric pair} $(\mathfrak{g},\mathfrak{h})$. Here, $\mathfrak{h}$ is the subalgebra of the fixed points of the differential of $\sigma$: 
  $$\mathfrak{h}=\mathfrak{g}^{\sigma}=\{X\in\mathfrak{g}\,|\,\sigma(X)=X\}.$$
  Note that throughout the paper we denote the differential of $\sigma$ by the same letter. 
  In the same way, we consider semisimple {\it 3-symmetric spaces}, which are defined by a similar condition (\ref{equo3}), but with the requirement that $\sigma$ is an automorphism of $G$ of order 3. These spaces were classified by Wolf and Gray \cite{wg1},\cite{wg2}. Various  interpretations of this class as well as its role in geometry are described in \cite{kow}.

Finally, let us explicitly formulate the basic assumptions on homogeneous spaces $G/H$ valid throughout this article.
\vskip6pt 
\noindent {\bf Basic assumptions on $G/H$.}  We consider homogeneous spaces $G/H$ of semisimple linear and connected real Lie groups $G$ and closed subgroups $H$ reductive in $G$, that is, $G/H$ are of reductive type.

Note that every semisimple symmetric space and every semisimple 3-symmetric space is of reductive type. This enables us to freely use all the results proved for this more general class. Although some lemmas and propositions in this paper may be formulated in a more general setting  the reader may keep in mind that they are valid   under the basic assumptions. 
     
\subsection{Hirzebruch-Kobayashi-Ono proportionality principle} \label{hirz}

In this subsection we assume that $G$ is a linear  semisimple Lie group, $H$ is a closed connected subgroup of $G$ and $\Gamma$ is a discrete subgroup of $G$ acting freely and properly on $G/H.$ Also let $m:=\dim G/H.$ Throughout this paper we use the relative cohomology $H^*(\mathfrak{g},\mathfrak{h},\mathbb{R})$ \cite{K1}. 

Recall that there is a natural map
$$\eta: (\Lambda(\mathfrak{g}/\mathfrak{h})^{\ast})^H\rightarrow \Omega(\Gamma\backslash G/H)$$
which induces a map in cohomology (denoted here by the same letter)
$$\eta: H^*(\mathfrak{g},\mathfrak{h},\mathbb{R})\rightarrow H^*(\Gamma\backslash G/H).$$
We don't describe it here referring to \cite{mor3}, Section 2 (and the references therein).

Assume that  $H$ is unimodular and that $\Gamma$ acts co-compactly on $G/H.$ It follows from the unimodularity of $H$ that
$$(\Lambda^{m}(\mathfrak{g}/\mathfrak{h})^{\ast})^{H}\neq 0.$$
Take any non-zero $\Phi \in (\Lambda^{m}(\mathfrak{g}/\mathfrak{h})^{\ast})^{H}$. We see that  $\eta(\Phi)$  is a volume form on $\Gamma \backslash G/H$. Since $\Gamma \backslash G/H$ is a compact manifold, $\eta([\Phi])\not=0$ in $H^{m}(\Gamma \backslash G/H,\mathbb{R})$. Therefore if $H$ is unimodular and $\Gamma$ acts co-compactly on $G/H$ then $\eta$ is injective in the top degree part.

Assume now that $G/H$ is of reductive type. Then $G/H$ has a dual $G_u/H_u$, which we now describe (see \cite{ko}, Section 3).  Let $G_{u}$ be a compact real form of a (connected) complexification $G^{\mathbb{C}}$ of $G$ and let $H_{u}$ be a compact real form of $H^{\mathbb{C}}\subset G^{\mathbb{C}}$ so that $H_{u}\subset G_{u}$. The space $G_{u}/H_{u}$ is called the  homogeneous space of compact type associated with $G/H$ (or dual to $G/H$). Groups $G_{u},$ $H_{u}$ are called the compact duals of $G$ and   $H,$ respectively. In such case
$$H^*(\mathfrak{g},\mathfrak{h},\mathbb{C})\cong H^*(G_{u}/H_{u},\mathbb{C})$$
so we can rewrite $\eta$ as a map from $H^*(G_{u}/H_{u},\mathbb{C})$ to $H^*(\Gamma\backslash G/H,\mathbb{C}).$
\begin{theorem}[\cite{ko}, Section 4, \cite{kkb}, Corollary 3.8] Let $G/H$ be a homogeneous space of a semisimple Lie group $G$ and of reductive type and assume that $H$ is connected. Assume that a discrete subgroup $\Gamma \subset G$ acts on $G/H$ freely and properly. Then there is a natural map
$$\eta\colon H^*(G_{u}/H_{u},\mathbb{C})\to H^*(\Gamma\backslash G/H,\mathbb{C}),$$
which sends the Euler class of $G_u/H_u$ to the Euler class of $\Gamma\backslash G/H$.
\label{ko}
\end{theorem}

\begin{remark} {\rm  Hirzebruch proved in \cite{H} the ``proportionality principle'' for the Chern numbers of compact Clifford-Klein forms of the bounded Hermitian symmetric domains and their compact duals.} 
\end{remark}

\subsection{Syndetic hull}\label{subsect:syndetic}
We will need the notion of a syndetic hull \cite{w-a}.
\begin{definition}
A {\it syndetic hull} of a closed subgroup $\Gamma$ of a Lie group $G$ is a closed subgroup $B$ of $G$ such that $B$ is connected, $B$ contains $\Gamma$ and $\Gamma\backslash B$ is compact.
\end{definition}
We will need the following results
\begin{theorem}[\cite{w-a}, Theorem 3.7]
Let $\Gamma$ be a discrete subgroup of a connected solvable linear Lie group $L$. If there is a compact subgroup $S$ of $\bar{\Gamma}$ and a compact subgroup $C$ of $L$ such that $SC$ is a maximal compact subgroup of the identity component $\bar{L}_{s}$ of $\bar{L}$, then some finite index subgroup $\Gamma_{0}$ of $\Gamma$ has a simply connected syndetic hull in $L$.
\label{thwit}
\end{theorem}

\begin{lemma}\label{lemma:syndetic} Let $G$ be a linear connected and semisimple Lie group. Let $\tilde{\Gamma}\subset G$ be a discrete solvable subgroup of $G.$ Then there exists a discrete solvable subgroup $\Gamma_{0}$ of finite index in $\tilde{\Gamma}$ with a simply connected syndetic hull $B\subset G$.
\end{lemma}
\begin{proof}  Since $G$ is linear, we can write $G=G_{\mathbb{R}}\subset GL(n,\mathbb{R})$. Take the identity component $L$ of the Zariski closure of $\tilde{\Gamma}$. Consider $\Gamma :=L\cap\tilde{\Gamma}$. We know that $\Gamma$ has a finite index in $\tilde{\Gamma}$. Apply Theorem \ref{thwit} to $L$ and $\Gamma$ (in this case one may take a maximal compact subgroup $C$ in $L$ and a trivial subgroup $S=\{e\}$ in $\bar{\Gamma}$). 
	\end{proof}

\section{Proof of Theorem \ref{twgg}}\label{sec:main}

\subsection{A Lie cohomology obstruction and proof of (i)}

In this subsection we find an obstruction to compact solvable Clifford-Klein forms using an idea very close to \cite{mor1}, Theorem 1.3. Assume that $G/H$ is of reductive type, $G$ is linear and semisimple. Assume that it admits a compact solvable Clifford-Klein form $\Gamma \backslash G/H$. By Lemma \ref{lemma:syndetic} we may assume without loss of generality  that $\Gamma$ has a simply connected syndetic hull $B$. Let $\mathfrak{g},$ $\mathfrak{h}$, $\mathfrak{b}$ be the Lie algebras of $G,H,B,$ respectively. In this subsection we will use the notion of the cohomological dimension $\operatorname{cd}_{\mathbb{R}}(\Gamma)$ of a discrete group $\Gamma$. We refer to \cite{ko-ma}. By definition
$$\operatorname{cd}_{\mathbb{R}}(\Gamma)=\sup\{n\in\mathbb{N}\,|\,H^n(\Gamma,\mathcal{A)}\not=0\,\text{for some left $\mathbb{R}[\Gamma]$-module}\,\mathcal{A}\}.$$ 
\begin{lemma}\label{lemma:man/b} Let $\Gamma\backslash G/H$ be a solvable compact Clifford-Klein form of a homogeneous space $G/H$ of reductive type. Let $\Gamma_H$ be a torsion-free co-compact lattice of $H$ and assume that $\Gamma$ has a simply connected syndetic hull $B$. Then $B\backslash G/\Gamma_H$ is a compact manifold of dimension $n=\dim\mathfrak{k}+\dim\mathfrak{p}_H$.
\end{lemma}
\begin{proof}
Let $G=\mathit{KAN}$ be the Iwasawa decomposition of $G$. By \cite{ko-ma}, Corollary 5.5 (b) the cohomological dimension of $\Gamma$ is given by the formula
$$\operatorname{cd}_{\mathbb{R}}\Gamma=\dim\mathfrak{p}-\dim\mathfrak{p}_H.$$
Recall that any simply connected solvable Lie group is diffeomorphic to $\mathbb{R}^{n}$ therefore $\operatorname{cd}_{\mathbb{R}}\Gamma=\dim B/\Gamma=\dim B$, which yields
$$\dim\,B\backslash G/\Gamma_H=\dim B\backslash G=\dim\mathfrak{g}-\dim\mathfrak{b}=$$
$$\dim\mathfrak{k}+\dim\mathfrak{p}-(\dim\mathfrak{p}-\dim\mathfrak{p}_H)=$$
$$\dim\mathfrak{k}+\dim\mathfrak{p}_H.$$
Also, $\dim\,B+\dim\mathfrak{p}_H=\dim\mathfrak{p}$. This equality implies the compactness of $B\backslash G/\Gamma_H$ by \cite{ow}, Theorem 3.4.
 \end{proof}

 \vskip6pt
 \centerline{\bf Proof of (i)}
 \vskip6pt 
Since $B$ admits a co-compact lattice it is unimodular. We have   shown in Section \ref{hirz} that  this implies
$H^{n}(\mathfrak{g},\mathfrak{b},\mathbb{R})\neq 0.$
Consider the following commutative diagram
\begin{center}
\begin{tikzpicture}
  \matrix (m) [matrix of math nodes,row sep=1em,column sep=2em,minimum width=2em]
  {
     H^{n}(\mathfrak{g}, \mathfrak{b},\mathbb{R}) & H^{n}(\mathfrak{g},\mathbb{R})  \\
     H^{n}(\Gamma_{H} \backslash G/B,\mathbb{R}) & H^{n}(\Gamma_{H} \backslash G) \\};
  \path[-stealth]
    (m-1-1) edge (m-2-1) edge node [left] {$\eta$} (m-2-1)
    (m-1-1.east|-m-1-2)        edge node [above] {$i^*$} (m-1-2)
    (m-2-1.east|-m-2-2) edge node [above] {$\pi^{\ast}$} (m-2-2)
    (m-1-2) edge  (m-2-2) edge node [right] {$\eta$} (m-2-2);

\end{tikzpicture}
\label{cd1} 
\end{center}
In this diagram $\pi^{\ast}$ is determined by a fibration with a contractible typical fiber
$$B\rightarrow \Gamma_{H} \backslash G \stackrel{\pi}{\rightarrow} \Gamma_{H} \backslash G/B,$$
which implies that $\pi^{\ast}$ is an isomorphism.
If $G/H$ admits a compact solvable Clifford-Klein form then the map
$$0\neq H^{n}(\mathfrak{g}, \mathfrak{b},\mathbb{R}) \stackrel{i^*}{\rightarrow} H^{n}(\mathfrak{g},\mathbb{R})$$
is injective.
Indeed, $\eta$ and $\pi^{\ast}$ are injective and it follows from the commutativity of the diagram above that   $i^*$ is injective as well. In particular one obtains $H^n(\mathfrak{g},\mathbb{R})\not=0$.  Note that  $n=\dim \Gamma_{H} \backslash G/B=\dim\mathfrak{k}+\dim\mathfrak{p}_H$ does not depend on the choice of $\Gamma$, by Lemma \ref{lemma:man/b}. The proof of (i) is complete.

\subsection{Proof of (ii)}
 The idea of the proof is similar to \cite{ko-ma}, Proposition 4.10. 
 Assume that $\operatorname{rank}\,G=\operatorname{rank}\,H$. 

Let $K$ and $K_H$ be as in Section \ref{sec:prelim}. Let $\Gamma\subset G$ be a discrete subgroup determining a compact Clifford-Klein form. By Selberg's lemma we may assume that $\Gamma$ is torsion-free. Recall that   $\Gamma$ acts co-compactly and freely on $G/H$. 
 Consider the fiber bundles
$$H/K_H\rightarrow G/K_H\rightarrow G/H,\,\text{and}\, K/K_H\rightarrow G/K_H\rightarrow G/K,$$
and
\begin{equation}
K/K_H\rightarrow \Gamma\backslash G/K_H\rightarrow \Gamma\backslash G/K.
\label{equo1} 
\end{equation}
Since $H/K_H$ is diffeomorphic to some euclidean space,  $\Gamma\backslash G/K_H$ and $\Gamma\backslash G/H$ have the same homotopy type.
Clearly, the classifying space for $\Gamma$ is  $B\Gamma=\Gamma\backslash G/K$, therefore $\chi(\Gamma)=\chi(\Gamma\backslash G/K).$
Looking at (\ref{equo1}) we write down
\begin{equation}
\chi(\Gamma\backslash G/H)=\chi(\Gamma\backslash G/K_H)=\chi(K/K_H)\cdot\chi(\Gamma). 
\label{equo2}
\end{equation}
 By assumption, $G$ and $H$ have equal ranks, hence the same is valid for $G_u$ and $H_u$. It is a classical result that $\chi(G_u/H_u)\not=0$ (see for example \cite{O} Theorem 2 on page 217) which implies that the Euler class of $T(G_u/H_u)$ is not zero by the Gauss-Bonnet-Chern theorem. By Theorem \ref{ko} the Euler class $e(T(\Gamma\backslash G/H))$ and the Euler characteristic $\chi(\Gamma\backslash G/H)$ also do not vanish.

On the other hand,  if $\Gamma$ is solvable then by Lemma \ref{lemma:syndetic} we may assume that $\Gamma$ has a simply connected syndetic hull $B$. As $B$ is a Lie group and $\Gamma$ is a discrete subgroup of $B$, the homogeneous space $B/ \Gamma$ is covered by $B$ and so the tangent bundle of $B/ \Gamma$ is trivial. Hence $\chi(B/\Gamma)=\chi(\Gamma)=0$ 
(for clear reasons we  assume that $\Gamma$ is infinite). 
Comparing this with (\ref{equo2}) yields a contradiction and so $\Gamma$ cannot be solvable.

\section{Proof of Theorem \ref{cwgg}}
To classify simple symmetric spaces that may admit compact solvable Clifford-Klein forms we take the Table A.1 in \cite{ok} and a corrected list from \cite{ok1} of such symmetric spaces that admit proper $SL(2,\mathbb{R})$ actions (if a non-compact semisimple symmetric space  admits a compact Clifford-Klein form it also admits a proper action of $SL(2,\mathbb{R})$, see \cite{ok}). Using Theorem \ref{twgg} (ii), we exclude rows that contain spaces of maximal rank and rows with spaces that do not admit compact Clifford-Klein forms listed in \cite{mor2}, Table 1. This way we obtain the following list of symmetric pairs: 
\begin{itemize}
\item exceptional pairs:
$$(\mathfrak{e}_{6(2)},\mathfrak{sp}(3,1)), \ (\mathfrak{e}_{6(-14)},\mathfrak{f}_{4(-20)}).$$
\item classical pairs:
$$(\mathfrak{so}(n,m),\mathfrak{so}(n-k,m)\times \mathfrak{so}(k)),\,\,\text{with}\,\,0<n\leq m;\,m,\, n \,\text{even}; \ k \ \text{odd},$$
$$(\mathfrak{su}(2p,2q),\mathfrak{sp}(p,q)),\,(\mathfrak{su}(2m-1,2m-1),\mathfrak{so}^*(4m-2)).$$
\end{itemize}
We need the following lemma.
\begin{lemma}\label{lemma:ex}
The non-compact symmetric space $G/H$ determined by the symmetric pair $(\mathfrak{e}_{6(-14)}, \mathfrak{f}_{4(-20)})$ does not admit compact solvable Clifford-Klein forms.
\end{lemma}
\begin{proof}
Assume that  a symmetric space $G/H$ admits a solvable compact Clifford-Klein form.  It follows from Theorem \ref{twgg} (i) that $H^{62}(\mathfrak{g})\not=0.$
We see that $H^{\ast}(\mathfrak{g}, \mathbb{R})\cong H^{\ast}(\mathfrak{g}_{u}, \mathbb{R}),$
 $\mathfrak{g}_{u}=\mathfrak{e}_{6},$ $\textrm{dim}\mathfrak{e}_{6}=78$. The  cohomology ring $H^*(\mathfrak{g}_u,\mathbb{R})$ has generators in degrees
$3, \ 9, \ 11, \ 15, \ 17, \ 23.$
It follows that $H^{62}(\mathfrak{g}_{u}, \mathbb{R})=H^{62}(\mathfrak{g}, \mathbb{R})=0.$
This is a contradiction. 
\end{proof}
The proof of Theorem \ref{cwgg} is complete.

\section{Proof of Theorem \ref{thm:3-sym}}
First we prove the following lemma
\begin{lemma}\label{lemma:ex2}
The 3-symmetric spaces $G/H$ determined by the 3-symmetric pairs
$$(\mathfrak{so}(4,4)/\mathfrak{su}(1,2)),\, (\mathfrak{so}(4,4)/\mathfrak{g}_{2(2)})$$
do not admit solvable compact Clifford-Klein forms.
\end{lemma}
\begin{proof}
Assume that $SO(4,4)/SU(1,2)$ admits a compact solvable Clifford-Klein form.
 It follows from Theorem \ref{twgg} (i) that 
$ H^{n}(\mathfrak{g}, \mathbb{R})\not=0$
for $n=16$. In the same way, if $SO(4,4)/G_{2(2)}$ admits a Clifford-Klein form, then $ H^{n}(\mathfrak{g}, \mathbb{R})\not=0$ for $n=20$. It is well known that 
$H^{\ast}(\mathfrak{g}, \mathbb{R})\cong H^{\ast}(\mathfrak{g}_{u}, \mathbb{R}),$
where $\mathfrak{g}_{u}$ is the compact dual of $\mathfrak{g}$. In our case $\mathfrak{g}_{u}=\mathfrak{so}(8),$ $\textrm{dim}\,\mathfrak{so}(8)=28$. The cohomology ring $H^*(\mathfrak{g}_u,\mathbb{R})$ has generators in degrees
$3, \ 7, \ 11, \ 7.$
It follows that $H^{16}(\mathfrak{g}, \mathbb{R})=0$ and $H^{20}(\mathfrak{g}, \mathbb{R})=0.$
We have arrived at a contradiction.
\end{proof}

The proof of Theorem \ref{thm:3-sym} is obtained by a combination of Lemma \ref{lemma:ex2}, Theorem \ref{twgg}  (ii) and the classification of non-compact 3-symmetric spaces of simple and absolutely simple real Lie groups in \cite{wg1},\cite{wg2} (Theorem 7.10 and Tables 7.11-7.14).

\end{document}